\documentclass[a4paper,reqno]{amsart}

\usepackage{amsmath}
\usepackage{amsfonts}
\usepackage{palatino}
\usepackage{a4wide}
\usepackage{amssymb}
\usepackage{amsthm}

\theoremstyle{plain}

\newtheorem{teo}{Theorem}[section]
\newtheorem{lemma}[teo]{Lemma}

\newtheorem{ackn}{Acknowledgments\!}

\theoremstyle{definition}

\theoremstyle{remark}

\newtheorem{rem}[teo]{Remark}

\numberwithin{equation}{section}

\def\R{{{\mathbb R}}}

\def\SSS{{{\mathcal S}}}

\def\NN{{{\mathbb N}}}

\def\trace{{\mathrm{tr}}}

\renewcommand{\a }{\alpha }
\renewcommand{\b }{\beta }
\renewcommand{\d}{\delta }
\newcommand{\D }{\Delta }

\newcommand{\g }{\gamma}

\newcommand{\G }{\Gamma }

\newcommand{\n }{\nabla }

\newcommand{\intbar}{\mathop{\int\makebox(-13.5,0){\rule[4pt]{.7em}{0.3pt}}%
\kern-6pt}\nolimits}
\newcommand{\wtilde }{\widetilde}

\newcommand{\be}{\begin{equation}}
\newcommand{\ee}{\end{equation}}

\newcommand{\ti}[1]{\wtilde{#1}}
\newcommand{\gt}{\ti{g}}
\newcommand{\p}{\partial}
\newcommand{\F}{\mathcal{F}}

\def\be{\begin{equation}}
\def\ee{\end{equation}}
\def\bea{\begin{eqnarray*}}
\def\bean{\begin{eqnarray}}
\def\eean{\end{eqnarray}}
\def\eea{\end{eqnarray*}}
\def\f{\frac}

\begin{document}

\title[On Perelman's Dilaton]
{On Perelman's Dilaton}

\author[Marco Caldarelli]{Marco Caldarelli}
\address[Marco Caldarelli]{}
\email[M. Caldarelli]{marco.caldarelli@gmail.com}

\author[Giovanni Catino]{Giovanni Catino}
\address[Giovanni Catino]{Dept. of Math., University of Pisa, Largo Bruno
Pontecorvo 5, Pisa, Italy, 56126}
\email[G. Catino]{catino@mail.dm.unipi.it}

\author[Zindine Djadli]{Zindine Djadli}
\address[Zindine Djadli]{Institut Fourier, 100 Rue des Maths, BP 74, St. Martin d'Heres, France, 38402}
\email[Z. Djadli]{Zindine.Djadli@ujf-grenoble.fr}

\author[Annibale Magni]{Annibale Magni}
\address[Annibale Magni]{SISSA -- International School for Advanced
  Studies, Via Beirut 2--4, 
Trieste, Italy, 34014}
\email[A. Magni]{magni@sissa.it}

\author[Carlo Mantegazza]{Carlo Mantegazza}
\address[Carlo Mantegazza]{Scuola Normale Superiore di Pisa, P.za Cavalieri 7, Pisa, Italy, 56126}
\email[C. Mantegazza]{mantegazza@sns.it}

\date{\today}

\begin{abstract} By means of a Kaluza--Klein type argument we show
  that the Perelman's ${\mathcal F}$--functional is the
  Einstein--Hilbert action in a space with extra ``phantom'' dimensions. 
In this way, we try to interpret some remarks of Perelman in the
introduction and at the end of the first section in his famous paper~\cite{perel1}.\\
As a consequence the Ricci flow (modified by a diffeomorphism and a time--dependent factor) is the
evolution of the ``real'' part of the metric under a 
constrained gradient flow of the Einstein--Hilbert gravitational action in higher dimension.
\end{abstract}

\maketitle
\tableofcontents

\section{Einstein--Hilbert Action and Perelman's ${\mathcal{F}}$--Functional}

Let $(M^m,g)$ and $(N^n,h)$ be two closed Riemannian manifolds of
dimension $m$ and $n$ respectively and let $f: M\rightarrow\R$ be a
smooth function on $M$. On the product manifold $\ti{M}=M\times N$ we
consider a metric $\ti{g}$ of the form $$\gt=e^{-Af}g\oplus
e^{-Bf}h,$$ where $A$ and $B$ are real constants. Notice that $\widetilde{g}$ is a conformal deformation of a warped product on $M$. We call the function $f$ {\em dilaton field}.

As a notation, we will use Latin indices, $i, j,\dots$ for the coordinates on $M$ (we will call them the "real" variables) and Greek indices, $\a, \b,\dots$, for the coordinates on $N$ (the "phantom" variables). Under these notations, clearly we have $\forall\, i,j\in\{1,\dots, m\}$ and $\forall\,\a,\b\in\{1,\dots, n\}$, 
$$
\gt_{i\a}=\gt^{i\a}=0\,,
$$ 
$$
\gt^{ij}=e^{Af}g^{ij}\,,\,\,\,\,\,\gt^{\a\b}=e^{Bf}h^{\a\b}\,.
$$
Let $\mu$, $\sigma$ and $\ti{\mu}$ be respectively the canonical volume measure on $M$, $N$ and $\ti{M}$. By definition of $\ti{g}$, it follows that $\ti{\mu}=e^{-\frac{Am+Bn}{2}}\mu\times\sigma$.

The Christoffel symbols of the metric $\gt$ are given by the formula
$$
\ti{\G}_{ab}^{c}=\f{1}{2}\gt^{cd}\left(\p_a\gt_{bd}+\p_b\gt_{ad}-\p_l\gt_{ab}\right)\,,
$$
where $a,b,\dots$ can be both real and phantom variables.

We have the following,
\begin{align*}
\ti{\G}_{ij}^k=&\,\f{1}{2}\gt^{kl}\left(\p_i\gt_{jl}+\p_j\gt_{il}-\p_l\gt_{ij}\right)\\
=&\,\f{1}{2}e^{Af}g^{kl}\left[e^{-Af}\left(\p_i g_{jl}+\p_j g_{il}-\p_l g_{ij}\right)-Ae^{-Af}\left(\p_i f g_{jl}+\p_j f g_{il}-\p_l f g_{ij}\right)\right]\\
=&\,\G_{ij}^k-\f{A}{2}\left(\p_i f \d^k_j+\p_j f \d^k_i-g^{kl} \partial_l f g_{ij}\right)\,.
\end{align*}
Using the fact that the metric $\gt$ is zero for a pair of ``mixed'' indexes and that the function $f$ depends only on the real 
variables, we get
\begin{align*}
\ti{\G}_{ij}^{\g}=&\,\frac{1}{2}\gt^{\g\b}\left(\p_i\gt_{j\b}+\p_j\gt_{i\b}-\p_{\b}\gt_{ij}\right)=0\,,\\
\ti{\G}_{\a i}^k=&\,\frac{1}{2}\gt^{kl}\left(\p_i\gt_{\a l}+\p_{\a}\gt_{il}-\p_l\gt_{i\a}\right)=0\,,\\
\ti{\G}_{i\beta}^{\g}=&\,\frac{1}{2}\gt^{\g\alpha}\left(\p_i\gt_{\a\b}+\p_\b\gt_{i\a}-\p_{\a}\gt_{i\b}\right)=
-\f{B}{2}\p_i f \d^{\g}_{\b}\,,\\
\ti{\G}_{\alpha\beta}^k=&\,\frac{1}{2}\gt^{kl}\left(\p_\a\gt_{l\b}+\p_{\b}\gt_{\a l}-\p_l\gt_{\a\b}\right)=
\f{B}{2}e^{(A-B)f}g^{kl}\partial_l f h_{\a\b}\,.
\end{align*}
Finally, a computation analogous to the one above gives $\ti{\G}_{\a\b}^{\g}=\G_{\a\b}^{\g}$.\\
Hence, summarizing
\begin{align}
\label{eq1}\ti{\G}_{ij}^k=&\,\G_{ij}^k-\f{A}{2}\left(\p_i f \d^k_j+\p_j f \d^k_i-g^{kl}\partial_l f g_{ij}\right)\\
\label{eq2}\ti{\G}_{ij}^{\a}=&\,\ti{\G}_{i\a}^k=0\\
\label{eq3}\ti{\G}_{\a\b}^k=&\,\f{B}{2}e^{(A-B)f}g^{kl}\partial_l f h_{\a\b}\\
\label{eq4}\ti{\G}_{i\b}^{\g}=&\,-\f{B}{2}\p_i f \d^{\g}_{\b}\\
\label{eq5}\ti{\G}_{\a\b}^{\g}=&\,\G_{\a\b}^{\g}\,.
\end{align}

We want now to compute the Ricci curvature of the metric $\gt$.\\
The Riemann tensor, as a $(1,3)$--tensor, is defined in terms of the derivatives of the Christoffel's symbols as follows 
$$
\ti{R}_{ab\,\,d}^{\,\,\,\,\,c}=\p_a\ti{\G}_{bd}^c-\p_b\ti{\G}_{ad}^c+\ti{\G}_{bd}^p\ti{\G}_{ap}^c-\ti{\G}_{ad}^p\ti{\G}_{bp}^c
$$
hence, the Ricci tensor is given by
$$
\ti{R}_{bd}=\p_a\ti{\G}_{bd}^a-\p_b\ti{\G}_{ad}^a+\ti{\G}_{bd}^p\ti{\G}_{ap}^a-\ti{\G}_{ad}^p\ti{\G}_{bp}^a\,.
$$
Using equations~\eqref{eq1}--~\eqref{eq5}, and computing in normal coordinates on both $M$ and $N$, we get the following 
\bea \ti{R}_{jl}&=&\p_i\ti{\G}_{jl}^i-\p_j\ti{\G}_{kl}^k-\p_j\ti{\G}_{\a l}^{\a}+\ti{\G}_{jl}^{k}\ti{\G}_{ki}^{i}+\ti{\G}_{jl}^{k}\ti{\G}_{\a k}^{\a}-\ti{\G}_{ij}^{k}\ti{\G}_{kl}^{i}-\ti{\G}_{\a j}^{\b}\ti{\G}_{\b l}^{\a}\\
&=&R_{jl}-\f{A}{2}\left(2\n^2_{jl}f-\D fg_{jl}\right)\\
&&+\f{Am}{2}\n^2_{jl}f+\f{Bn}{2}\n^2_{jl}f\\
&&+\f{A^2m}{4}\left(2df_j df_l-|\n f|^2 g_{jl}\right)+\f{ABn}{4}\left(2 df_jdf_l -|\n f|^2 g_{jl}\right)\\
&&-\f{A^2}{4}\left[(m+2)df_j df_l -2|\n f|^2 g_{jl}\right]-\f{B^2n}{4} df_j df_l\,,
\eea
that is, collecting similar terms,
\bean\label{ricci} \ti{R}_{jl}&=& R_{jl}+\n^2_{jl} f \left(\f{Am+Bn}{2}-A\right)
+\f{A}{2}g_{jl}\left[\D f-|\n f|^2\left(\f{Am+Bn}{2}-A\right)\right]\\
&&+\f{1}{4} df_j df_l \left(2ABn+(m-2)A^2-B^2n\right)\,.\nonumber
\eean
On the other hand, for the phantom indexes, we get
\bea \ti{R}_{\b\g}&=&\p_{\a}\ti{\G}_{\b\g}^{\a}-\p_{\g}\ti{\G}_{\a\b}^{\a}+\p_k\ti{\G}_{\b\g}^k+\ti{\G}_{\b\g}^{k}\ti{\G}_{\a\g}^{\a}+\ti{\G}_{\b\g}^{k}\ti{\G}_{ik}^{i}-\ti{\G}_{\a\g}^{k}\ti{\G}_{\b k}^{\a}-\ti{\G}_{k\g}^{\a}\ti{\G}_{\a\b}^{k}\\
&=& R_{\b\g}+\f{B}{2}e^{(A-B)f}h_{\b\g}\left(\D f + (A-B)|\n f|^2\right)\\
&&-\f{B^2 n}{4}e^{(A-B)f}h_{\b\g}|\n f|^2-\f{ABm}{4}e^{(A-B)f}h_{\b\g}|\n f|^2\\
&&+\f{B^2}{4}e^{(A-B)f}h_{\b\g}|\n f|^2+\f{B^2}{4}e^{(A-B)f}h_{\b\g}|\n f|^2\,,
\eea
that is,
\begin{equation}
\ti{R}_{\b\g}=R_{\b\g}+\f{B}{2}e^{(A-B)f}h_{\b\g}\left[\D f-|\n
  f|^2\left(\f{Am+Bn}{2}-A\right)\right]\,.\label{riccih}
\end{equation}
Finally, it is easy to see that the mixed terms of the Ricci tensor of $\gt$ vanish, $\ti{R}_{i\a}=0$.

From this computation we get then the formula for the scalar curvature of $\gt$, 
\bea \ti{R}&=&e^{Af}R^{M}+e^{Bf}R^{N}+e^{Af}\D f(Am+Bn-A)\\
&&+\f{e^{Af}}{4}|\n f|^2\left(4ABn-2ABmn+3mA^2-2A^2-m^2A^2-B^2n-B^2n^2\right)\,.
\eea 
where $R^M$ and $R^N$ are respectively the scalar curvatures of $(M,g)$ and $(N,h)$.

We make now the following ansatz:
\begin{equation*}
\label{Constr1}\tag{C1} 2ABn+(m-2)A^2-B^2n=0
\end{equation*}
and
\begin{equation*}
\label{Constr2}\tag{C2}\frac{Am+Bn}{2}-A=1 \qquad \Longleftrightarrow\qquad A(m-2)+Bn=2\,.
\end{equation*}

\begin{rem}\label{motivatic} 
We spend some words to motivate our choice of the constants $A$ and $B$, which
we guess it is not very clear at this point.\\
Condition (C1) is assumed in order to make vanish from the expression
of $\widetilde{R}_{ij}$ the term in $df\otimes df$ that otherwise 
appears in doing the flow by the gradient of the functional
$\int_{\widetilde{M}}\widetilde{R}\,d\widetilde{\mu}$ (see
Section~\ref{other} and Remark~\ref{remlambda}).\\
The second condition, that clearly also simplifies both
$\widetilde{R}_{ij}$ and $\widetilde{R}_{\alpha\beta}$, is instead more related to Perelman's
${\mathcal F}$--functional. In writing the functional
$\int_{\widetilde{M}}\widetilde{R}\,d\widetilde{\mu}$ as an integral
on $M$ with respect to the measure $\mu$ we will see that the only way
to get the factor $e^{-f}$ is to assume condition (C2).
\end{rem}

\begin{lemma}
If $m+n>2$, we can always find two non zero constants $A$ and $B$ satisfying these two
conditions.
\end{lemma}
\begin{proof} Notice that $A=0$ implies $B=0$. If $B\not=0$, dividing
  both sides of condition (C1) by $B^2$, it can be expressed in the
  following form for $\theta=A/B$,
\begin{equation*}
\label{Constr1*}\tag{C1$^*$}(m-2)\theta^2+2n\theta-n=0\,.
\end{equation*}
If $m\not=2$, this second degree equation for $\theta$ has always two 
solutions for every choice of the dimensions $m,n\in\NN$, which would
coincide only in the case $m=n=1$, that we excluded.\\
Notice also that the two solutions have opposite signs. Precisely, they are
$$
\theta=\frac{-n\pm\sqrt{n(n+m-2)}}{m-2}
$$
and in the special case $n=1$, we have $\theta=\frac{-1\pm\sqrt{m-1}}{m-2}$.\\
If $m=2$ we have only one solution of equation (C1$^*$) which is $\theta=1/2$.\\
Then, condition (C2) is equivalent to $\theta(m-2)+n=2/B$ which can be
fulfilled, by homogeneity, if $\theta(m-2)+n\not=0$. If this happen, we would have
$$
0=\theta^2(m-2)+2n\theta-n=n\theta-n
$$
which would imply $\theta=1$. Hence, $m-2+n=0$ and $m=n=1$.
\end{proof}

Under assumptions (C1) and (C2), the last term of $\ti{R}_{jl}$ in
formula~\eqref{ricci} cancels out and many coefficients becomes one. We 
get indeed the following ``smooth`` formulas for the components of the Ricci tensor of $\gt$,
\be\label{Rics}\ti{R}_{jl}=R_{jl}+\n^2_{jl}f+\f{A}{2}g_{jl}\left(\D f-|\n f|^2\right),\ee
\be\label{Rict}\ti{R}_{\b\g}=R_{\b\g}+\f{B}{2}e^{(A-B)f}h_{\b\g}\left(\D f-|\n f|^2\right)\,.\ee
Then, the scalar curvature of $\gt$ becomes 
\be\label{Scal}\ti{R}=e^{Af}R^M+e^{Bf}R^N+e^{Af}\left(\D f (A+2)-|\n f|^2(A+1)\right)\,.\ee

From this last formula, it follows immediately the relation between
the Einstein--Hilbert action 
functional $\SSS$ on $\ti{M}$ 
and the Perelman's ${\mathcal{F}}$--functional, see~\cite{perel1},
$$
\F(g,f)=\int_M (R^M+\vert\nabla f\vert^2)e^{-f}\,d\mu\,.
$$

\begin{teo}\label{Hilbert} Let $(M^m,g)$ and $(N^n,h)$ be two closed Riemannian manifolds of dimension $m$ and $n$ respectively, with $m+n>2$ and let $f: M\rightarrow\R$ be a smooth function on $M$. 
On the product manifold $\ti{M}=M\times N$ consider the metric $\ti{g}$ of the form 
$$
\gt=e^{-Af}g\oplus e^{-Bf}h\,,
$$
where $A$ and $B$ are constants satisfying conditions (\ref{Constr1}) and (\ref{Constr2}).\\
Then the following formula holds 
\be\SSS(\gt)=\int_{\ti{M}}\ti{R}\,d\ti{\mu}={\mathrm{Vol}}(N,h)\F(g,f)
+\left(\int_{M}e^{(B-A-1)f}\,d\mu\right)\int_{N}R^N\,d\sigma\ee
In particular, if $(N,h)$ has zero total scalar curvature and unit volume, we get $\SSS(\gt)=\F(g,f)$.
\end{teo}
\begin{proof}
We simply compute
\begin{align*}
\int_{\ti{M}}\ti{R}\,d\ti{\mu}=&\,\int_M\int_N e^{-\frac{Am+Bn}{2}f}\ti{R}\,d\mu\,d\sigma\\
=&\,\int_M\int_N e^{-(1+A)f}\left[e^{Af}R^M+e^{Bf}R^N+e^{Af}\left(\D f (A+2)-|\n f|^2(A+1)\right)\right]\,d\mu\,d\sigma\\
=&\,\int_M\int_N e^{-(1+A)f}e^{Bf}R^N\,d\mu\,d\sigma\\
&\,+\int_M\int_N\left[R^M+\D f (A+2)-|\n f|^2(A+1)\right]e^{-f}\,d\mu\,d\sigma\\
=&\,\left(\int_{M}e^{(B-A-1)f}\,d\mu\right)\int_{N}R^N\,d\sigma\\
&\,+{\mathrm{Vol}}(N,h)\int_M\left[R^M+\D f (A+2)-|\n f|^2(A+1)\right]e^{-f}\,d\mu\\
=&\,\left(\int_{M}e^{(B-A-1)f}\,d\mu\right)\int_{N}R^N\,d\sigma\\
&\,+{\mathrm{Vol}}(N,h)\int_M\left(R^M+|\n f|^2\right)e^{-f}\,d\mu
\end{align*}
where in the last passage we integrated by parts the Laplacian term.
\end{proof}

\section{The Associated Flow}

Under assumptions (C1) and (C2), we have
\be\ti{R}_{jl}=R_{jl}+\n^2_{jl}f+\f{A}{2}g_{jl}\left(\D f-|\n f|^2\right)\,,\qquad\ti{R}_{i\alpha}=0\,,\ee
\be\ti{R}_{\b\g}=R_{\b\g}+\f{B}{2}e^{(A-B)f}h_{\b\g}\left(\D f-|\n f|^2\right)\ee
and
\be\ti{R}=e^{Af}R^M+e^{Bf}R^N+e^{Af}\left(\D f (A+2)-|\n f|^2(A+1)\right)\,.\ee

Suppose we have a manifold $\ti{M}=M\times N$ with a time dependent metric 
$\ti{g}(t)$ for $t\in[0,T]$.\\
If the initial metric is a warped product $\ti{g}=\widehat{g}\oplus \varphi h$ with $\varphi:M\to\R$ a smooth function, $(N,h)$ Ricci--flat and of unit volume, we consider the motion by the gradient of the Einstein--Hilbert action with the constraint that the  measure $\varphi^{-\theta}\ti{\mu}$ is fixed, where $\theta$ comes from condition (C1$^*$) and $A$, $B$ are the relative constants satisfying conditions (C1) and (C2) above.\\
Suppose there exists a unique solution of this flow, preserving the
warped product. 
We can assume that for every $t\in[0,T]$ we have $\ti{g}(t)=\widehat{g}(t)\oplus \varphi(t) h(t)$ with $(N,h(t))$ always of volume 1.\\
Writing down the evolution of $h$ we see that it moves only by multiplication by a positive factor, as we assumed that $(N,h(t))$ is of unit volume, we then conclude that the metric $h(t)$ has to be constant equal to the initial $h$.
Setting $f=-\frac{1}{B}\log{\varphi}$ which implies $\varphi=e^{-Bf}$ and $\varphi^{-\theta}=e^{-Af}$, and we can write
$\ti{g}=e^{-Af}g\oplus e^{-Bf} h$ where $g(t)=e^{Af}\widehat{g}(t)$. Clearly, also $\ti{g}=\varphi^\theta g\oplus \varphi h$.

Denote with $\d\gt$, $\d g$ and $\d f$ the variations of $\gt$, $g$ and $f$
respectively. Then we have,
$$
\d\gt=e^{-Af}\left(\d g-Ag\d f\right)\oplus
e^{-Bf}\left(-Bh\d f\right)\,.
$$
and the constraint, in terms of these variations becomes $\d f=\trace_g\delta
g/2$. Keeping in mind that $(N,h)$ is Ricci--flat, we get
\begin{align*}
\delta\int_{\ti{M}}2\ti{R}\,d\ti{\mu}
=&\int_{\ti{M}}\langle-2\ti{Ric}+\ti{R}\gt\,|\,\d\gt\rangle\,d\ti{\mu}\\
=&\int_{\ti{M}}\langle-2\ti{Ric}+\ti{R}\gt\,|\,e^{-Af}\left(\d g-Ag\d f\right)\oplus e^{-Bf}\left(-Bh\d f\right)\rangle\,d\ti{\mu}\\
=&\int_{\ti{M}}\langle-2(Ric^M+\n^2 f)\,|\,\d g\rangle e^{-Af}\,d\ti{\mu}\\
&\,+\int_{\ti{M}}\left[-A(\D f-|\n f|^2)+(R^M+\D f (A+2)-|\n f|^2
(A+1))\right]\trace_g(\d g) e^{-Af}\,d\ti{\mu}\\
&\,+\int_{\ti{M}}\Bigl\langle-2\Bigl[Ric^M+\n^2 f+\f{A}{2}g(\D f -|\n
f|^2)\Bigr]\,\Bigl|\,-Ag\d f\Bigr\rangle e^{-Af}\,d\ti{\mu}\\
&\,+\int_{\ti{M}}\left[R^M+\D f (A+2)-|\n f|^2 (A+1)\right](-Am\d
f) e^{-Af}\,d\ti{\mu}\\
&\,+\int_{\ti{M}}\left[-B(\D f - |\n f|^2)\right](-Bn\d
f) e^{-Af}\,d\ti{\mu}\\
&\,+\int_{\ti{M}}\left[R^M+\D f (A+2)-|\n f|^2 (A+1)\right](-Bn\d
f) e^{-Af}\,d\ti{\mu}\\
=&\int_{\ti{M}}\langle-2(Ric^M+\n^2 f)\,|\,\d g\rangle e^{-Af}\,d\ti{\mu}\\
&\,-\f{1}{2}\int_{\ti{M}}(\D f -|\n f|^2)(ABn+2A-B^2n)\trace_g(\d
g)e^{-Af}\,d\ti{\mu}\\
=&\int_{\ti{M}}\langle-2(Ric^M+\n^2 f)\,|\,\d g\rangle e^{-Af}\,d\ti{\mu}\\
=&-2\int_{{M}}\langle Ric^M+\n^2 f\,|\,\d g\rangle e^{-f}\,d\mu\,,
\end{align*}
since, by conditions (C1) and (C2), it follows $ABn+2A-B^2n=0$.\\
Hence, the system
$$
\begin{cases}
\delta g=-2(Ric^M+\n^2 f)\\
\delta f=-\Delta f -R^M
\end{cases}
$$
represents the constrained gradient of the
Einstein--Hilbert action functional. The associated flow of the metric $\widetilde{g}=e^{-Af}g\oplus e^{-Bf}h$ is
described by
$$
\begin{cases}
\partial_t g=-2(Ric^M+\n^2 f)\\
\partial_t h=0\\
\partial_t f=-\Delta f -R^M\,,
\end{cases}
$$
that is, $g$ evolves by the ``modified'' Ricci flow.\\
Following Perelman~\cite{perel1}, modifying the pair $(g,f)$ by a
suitable diffeomorphism,
we get a solution of
$$
\begin{cases}
\partial_t g=-2Ric^M\\
\partial_t f=-\Delta f +\vert\nabla f\vert^2 - R^M
\end{cases}
$$
hence, up to a factor and a diffeomorphism, the spatial part of the
metric $\widetilde{g}$ moves according to the Ricci flow 
($g$ is equal to the spatial part of $\widetilde{g}$ times the factor $e^{Af}$).

\section{Other Flows}\label{other}

It is interesting to see what functionals and flows one can get by varying the
constants $A$ and $B$.\\
Supposing that $(N,h)$ has unit volume and zero total scalar curvature, we computed,
\bea\ti{R}_{jl}&=& R_{jl}+\n^2_{jl} f \left(\f{Am+Bn}{2}-A\right)
+\f{A}{2}g_{jl}\left[\D f-|\n f|^2\left(\f{Am+Bn}{2}-A\right)\right]\\
&&+\f{1}{4} df_j df_l \left(2ABn+(m-2)A^2-B^2n\right)\\
\ti{R}_{\b\g}&=&R_{\b\g}+\f{B}{2}e^{(A-B)f}h_{\b\g}\left[\D f-|\n
  f|^2\left(\f{Am+Bn}{2}-A\right)\right]
\eea 
Assuming the condition $\f{Am+Bn}{2}-A=1$ we have
\bea\ti{R}_{jl}&=& R_{jl}+\n^2_{jl} f +\f{A}{2}g_{jl}\left[\D f-|\n f|^2\right]\\
&&+\f{1}{4} df_j df_l \left(2ABn+(m-2)A^2-B^2n\right)\\
\ti{R}_{\b\g}&=&R_{\b\g}+\f{B}{2}e^{(A-B)f}h_{\b\g}\left[\D f-|\n
  f|^2\right]\\
\ti{R}&=&e^{Af}R^{M}+e^{Bf}R^{N}+e^{Af}\Delta f\\
&&+e^{Af}\left(\frac{Am+Bn}{2}\right)(\Delta
f-\vert\nabla
f\vert^2)+\f{e^{Af}}{4}\left(2ABn+(m-2)A^2-B^2n\right)\vert \nabla f\vert^2\\
&=&e^{Af}R^{M}+e^{Bf}R^{N}+e^{Af}\Delta f\\
&&+e^{Af}(A+1)(\Delta f-\vert\nabla
f\vert^2)+\f{e^{Af}}{4}\left(2ABn+(m-2)A^2-B^2n\right)\vert \nabla f\vert^2\,.
\eea 
Hence,
\begin{align*}
\int_{\ti{M}}\ti{R}\,d\ti{\mu}=&\,\int_M\int_N e^{-\frac{Am+Bn}{2}f}\ti{R}\,d\mu\,d\sigma\\
=&\,\int_M\int_N e^{-(1+A)f}\left[e^{Af}R^{M}+e^{Bf}R^{N}+e^{Af}\Delta
  f\right]\,d\mu\,d\sigma\\
&\,+\int_M\int_N e^{-(1+A)f} e^{Af}(A+1)(\Delta
f-\vert\nabla f\vert^2)\,d\mu\,d\sigma\\
&\,+\int_M\int_N e^{-(1+A)f} \f{e^{Af}}{4}\left(2ABn+(m-2)A^2-B^2n\right)\vert \nabla f\vert^2\,d\mu\,d\sigma\\
=&\,\int_M\int_N e^{-(1+A)f}e^{Bf}R^N\,d\mu\,d\sigma\\
&\,+\int_M\int_N\left[R^M+\Delta f+\frac{1}{4}|\n f|^2(2ABn+(m-2)A^2-B^2n)\right]e^{-f}\,d\mu\,d\sigma\\
=&\,\int_M\left[R^M+\vert\nabla f\vert^2+\frac{1}{4}|\n
  f|^2(2ABn+(m-2)A^2-B^2n)\right]e^{-f}\,d\mu\\
=&\,\,{\mathcal{F}}(g,f)+Z_{m,n}(A,B)\int_M\vert\nabla
f\vert^2e^{-f}\,d\mu\,,
\end{align*}
with $Z_{m,n}(A,B)=(2ABn+(m-2)A^2-B^2n)/4$.\\

We want to see what are the possible values of $Z_{m,n}$, we recall
that we have the constraint $A(m-2)+Bn=2$.\\
We change variables as $x=A$ and $y=(B-A)$ so the constraint becomes
$(m+n-2)x+ny=2$ and $4Z_{m,n}(A,B)=(m+n-2)x^2-ny^2$. As
$y=[2-x(m+n-2)]/n$ we get (like before we assume $m+n>2$),
\begin{align*}
4Z_{m,n}(A,B)=&\,(m+n-2)x^2-n\left(\frac{2-x(m+n-2)}{n}\right)^2\\
=&\,(m+n-2)x^2-(4+x^2(m+n-2)^2-4x(m+n-2))/n\\
=&\,x^2[(m+n-2)-(m+n-2)^2/n]+4x(m+n-2)/n-4/n\\
=&\,-x^2\frac{(m+n-2)(m-2)}{n}+x\frac{4(m+n-2)}{n}-\frac{4}{n}\,.
\end{align*}
In the special case $m=2$, we have $B=2/n$ and $A$ ``free'', then
\begin{align*}
Z_{m,n}(A,B)=\frac{x(m+n-2)-1}{n}=\frac{An-1}{n}=A-1/n
\end{align*}
which can take every real value as $x$ can vary from $-\infty$ to
$+\infty$.\\
If instead, $m>2$ the expression
\begin{equation*}
Z_{m,n}(A,B)=-A^2\frac{(m+n-2)(m-2)}{4n}+A\frac{m+n-2}{n}-\frac{1}{n}\,.
\end{equation*}
is a second degree polynomial in $A\in\R$ with negative leading
coefficient, so it can vary only between $-\infty$ and some
maximum. By a straightforward computation one sees that such a maximum is
given by $1/(m-2)$, which is independent of the dimension $n$.\\
This means that by a suitable choice of the constants $A$ and $B$ one
has
$$
\SSS(\gt)=\int_{\ti{M}}\ti{R}\,d\ti{\mu}=\int_M(R^M+(\lambda+1)\vert\nabla
f\vert^2)e^{-f}\,d\mu\,,
$$
for every $\lambda\in\left(-\infty,\frac{1}{m-2}\right]$. Notice that (if $m>2$), with the exception of $\lambda=1/(m-2)$ one has always two possible choices of pairs of constants $(A,B)$ for every value $\lambda$. 

When $\lambda\not=0$ as 
$$
\ti{R}_{jl}=R_{jl}+\n^2_{jl} f +\f{A}{2}(\D f-|\n
  f|^2)g_{jl}+\lambda(df\otimes df)_{jl}\,,
$$
the associated flow is substantially different from 
the (modified) Ricci flow, indeed if as before 
$\d f =\f{1}{2}\trace_g (\d g)$ and $(N,h)$ is Ricci--flat, we get
\begin{align*}
\d\int_{\ti{M}}2\ti{R}\,d\ti{\mu}
=&\int_{\ti{M}}\langle-2\ti{Ric}+\ti{R}\gt\,|\,\d\gt\rangle\,d\ti{\mu}\\
=&\int_{\ti{M}}\langle-2\ti{Ric}+\ti{R}\gt\,|\,e^{-Af}\left(\d g-Ag\d f\right)\oplus e^{-Bf}\left(-Bh\d f\right)\rangle\,d\ti{\mu}\\
=&\int_{\ti{M}}\langle-2(Ric^M+\n^2 f+\lambda df\otimes df)\,|\,\d g\rangle e^{-Af}\,d\ti{\mu}\\
&\,+\int_{\ti{M}}\left[-A(\D f-|\n f|^2)+(R^M+\D f (A+2)-|\n f|^2
(A+1))\right]\trace_g(\d g) e^{-Af}\,d\ti{\mu}\\
&\,+\int_{\ti{M}}\Bigl\langle-2\Bigl[Ric^M+\n^2 f+\f{A}{2}g(\D f -|\n
f|^2)\Bigr]\,\Bigl|\,-Ag\d f\Bigr\rangle e^{-Af}\,d\ti{\mu}\\
&\,+\int_{\ti{M}}\left[R^M+\D f (A+2)-|\n f|^2 (A+1)\right](-Am\d
f) e^{-Af}\,d\ti{\mu}\\
&\,+\int_{\ti{M}}\left[-B(\D f - |\n f|^2)\right](-Bn\d
f) e^{-Af}\,d\ti{\mu}\\
&\,+\int_{\ti{M}}\left[R^M+\D f (A+2)-|\n f|^2 (A+1)\right](-Bn\d
f) e^{-Af}\,d\ti{\mu}\\
=&-2\int_{{M}}\langle Ric^M+\n^2 f+\lambda df\otimes df\,|\,\d g\rangle e^{-f}\,d\mu\,.
\end{align*}
Hence, as before, the system
$$
\begin{cases}
\delta g=-2(Ric^M+\n^2 f+\lambda df\otimes df)\\
\delta f=-\Delta f -R^M  -\lambda\vert \nabla f\vert^2 
\end{cases}
$$
represents the constrained gradient of the
Einstein--Hilbert action functional and the associated flow is
$$
\begin{cases}
\partial_t g=-2(Ric^M+\n^2 f+\lambda df\otimes df)\\
\partial_t f=-\Delta f -R^M-\lambda\vert \nabla f\vert^2\,.
\end{cases}
$$
\begin{rem}\label{remlambda} Like in the Ricci flow, the flow $\partial_t g=-2(Ric+\n^2
  f+\lambda df\otimes df)$ can be modified by a diffeomorphism to the
  flow $\partial_t g=-2(Ric+\lambda df\otimes df)$. The extra term
  $df\otimes df$ instead, cannot be ``canceled'' in this way as $\nabla^2f$.
\end{rem}

Notice that, as in Perelman's work, immediately one gets the
monotonicity of the relative ${\mathcal{F}}$--functional along this
flow.
$$
\frac{d\,}{dt}\int_M(R^M+(\lambda+1)\vert\nabla
f\vert^2)e^{-f}\,d\mu=-2\int_{{M}}\vert Ric^M+\n^2 f+\lambda df\otimes
df\vert^2 e^{-f}\,d\mu\,.
$$
\begin{rem} Compare the material of this section with the Ph.D~Thesis
  of List~\cite{listphd}.
\end{rem}

\begin{ackn} The authors are grateful to the CRM of Barcelona for the warm
  hospitality during the period this research was carried on.\\
Z. Djadli is supported by the project ANG "Flots et Op\'erateurs
G\'eom\'etriques ANR-07-BLAN-0251-01".\\
A. Magni is partially supported by the ESF Programme 
"Methods of Integrable Systems, Geometry, Applied
Mathematics" (MISGAM) and Marie Curie RTN "European Network in Geometry,
Mathematical Physics and Applications" (ENIGMA).
\end{ackn}

\bibliographystyle{amsplain}
\bibliography{dilaton}

\providecommand{\bysame}{\leavevmode\hbox to3em{\hrulefill}\thinspace}
\providecommand{\MR}{\relax\ifhmode\unskip\space\fi MR }
\providecommand{\MRhref}[2]{%
  \href{http://www.ams.org/mathscinet-getitem?mr=#1}{#2}
}
\providecommand{\href}[2]{#2}
\begin{thebibliography}{1}

\bibitem{listphd}
B.~List, \emph{Evolution of an extended {R}icci flow system}, Ph.D. thesis,
  Max--Planck--Instituts fur Gravitationsphysik (Albert Einstein Institut),
  Potsdam, 2005.

\bibitem{perel1}
G.~Perelman, \emph{The entropy formula for the ricci flow and its geometric
  applications}, ArXiv Preprint Server -- http://arxiv.org, 2002.

\end{thebibliography}

\end{document}